\documentclass[12pt,reqno]{amsart}
\usepackage{amssymb}
\usepackage{graphicx}
\usepackage{mathrsfs}
\usepackage{amsfonts}
\usepackage{caption}
\usepackage{textcomp}
\usepackage{amsfonts,amssymb}
\newtheorem{thm}{Theorem}
\newtheorem{lemma}{Lemma}
\newtheorem*{remark}{Remark}

\begin{document}

\noindent
{\Large {\bf Iterative solution of a nonlinear static beam equation  }}\\[2mm]
{\large {\bf Givi Berikelashvili$^{1, 2}$, Archil Papukashvili$^{3, 4}$, Giorgi Papukashvili$^{1,5}$, and Jemal Peradze$^{1,3}$}}\\
{\sl $^1$ Department of Mathematics of Georgian Technical University, Tbilisi, Georgia}\\
{\sl ${^2}$ A.Razmadze Mathematical Institute,  I. Javakhishvili Tbilisi State University, Tbilisi, Georgia}\\
{\sl $^3$ Faculty of Exact and Natural Sciences, I. Javakhishvili Tbilisi State University, Tbilisi, Georgia} \\
{\sl $^4$ I. Vekua Institute of Applied Mathematics, I. Javakhishvili Tbilisi State University, Tbilisi, Georgia}\\
{\sl $^5$ V. Komarovi Tbilisi Physics and Mathematics \textnumero 199 Public School, Tbilisi, Georgia}
\vspace{0.5cm}

\textbf{ Abstract}
The paper deals with a boundary value problem for the nonlinear integro-differential equation $u^{\prime\prime\prime\prime}-m\left(\int_0^l {u^\prime}^2dx\right)u^{\prime\prime}=f(x,u,u^\prime), \; m(z)\geq \alpha>0, \; 0\leq z <\infty$, modelling the static state of the Kirchhoff beam. The problem is reduced to a nonlinear integral equation which is solved using the Picard iteration method. The convergence of the iteration process is established and the error estimate is obtained.
\noindent

\noindent


\section{Statement of the Problem}

Let us consider the nonlinear beam equation
\begin{equation}\label{eq:bpp-01}
  u^{\prime\prime\prime\prime}(x)-m\left(\int_0^l {u^\prime}^2(x)dx\right)u^{\prime\prime}(x)=f(x,u(x),u^\prime(x)), \quad x\in \left( 0,l \right),
\end{equation}
with the conditions
\begin{equation}\label{eq:bpp-02}
  u(0)=u(l)=0,\; u^{\prime\prime}(0)=u^{\prime\prime}(l)=0.
\end{equation}
Here $u=u(x)$ is the displacement function of length $l$ of the beam subjected to the action of a force given by the function $f(x,u, u^\prime)$, the function $m(z)$,
\begin{equation}\label{eq:bpp-03}
  m(z) \geq \alpha>0, \quad 0\leq z<\infty,
\end{equation}
describes the type of a relation  between stress and strain. Namely, if the function $m(z)$ is linear, this means that this relation is consistent with Hooke's linear law, while otherwise we deal with material nonlinearities.

Equation \eqref{eq:bpp-01} is the stationary problem associated with the equation
$$
u_{tt}+u_{xxxx}-m\left(\int_0^l u^2_x dx\right)u_{xx}=f(x,t,u,u_x),
$$
$$
m(z)\geq \text{const} >0,
$$
which for the case where $m(z)=m_0+m_1z,\; m_0,m_1 >0,$ and $f(x,t,u,u_x)=0,$ was proposed by Woinowsky-Krieger [11] as a model of deflection of an extensible dynamic beam with hinged ends. The nonlinear term $\int_0^lu^2_xdx$  was for first time used by Kirchhoff [3] who generalized D'Alembert's classical linear model. Therefore \eqref{eq:bpp-01} is frequently called a Kirchhoff type equation for a static beam.

The problem of construction of numerical algorithms and estimation of their accuracy for equations of type \eqref{eq:bpp-01} is investigated in [1], [5], [8] and [9].
In [4], the existence of a solution of problem \eqref{eq:bpp-01},\,\eqref{eq:bpp-02} is proved when the right-hand part of equation is written in the form $q(x)f(x,u,u')$, where $f\in C([0,l]\times [0,\infty)\times \mathbb{R})$ is a nonnegative function and $q\in C[0,l]$ is a positive function.

In the present paper, in order to obtain an approximate solution of the problem \eqref{eq:bpp-01},\eqref{eq:bpp-02},  an approach is used, which differs from those applied in the above-mentioned references. It consists in reducing the problem \eqref{eq:bpp-01},\eqref{eq:bpp-02} by means of Green's function to a nonlinear integral equation, to solve which we use the iterative process. The condition for the convergence of the method is established and its accuracy is estimated.

The Green's function method with a further iteration procedure has been applied by us previously also to a nonlinear problem for the axially symmetric Timoshenko plate [6].

\vskip5mm
\noindent
\section{ Assumptions}
\vskip 2mm
Let us assume that besides \eqref{eq:bpp-03} the function $m(z)$ also satisfies the Lipschitz condition
$$
|m(z_1)-m(z_2)| \leq l_1|z_2-z_1|,\quad 0\leq z_1, z_2 <\infty,\quad l_1=\text{const} >0.
$$
Suppose that $f(x,u,v)\in L_2\left((0, l), \mathbb{R},\mathbb{R}\right)$ and, additionally, that the inequalities
\begin{equation}\label{eq:bpp-04}
  |f(x,u,v)|\leq \sigma_1(x)+\sigma_2(x)\,|u|+\sigma_3(x) \, |v|, \,
\end{equation}
\begin{equation}\label{eq:bpp-05}
  |f(x,u_2,v_2)-f(x,u_1,v_1)|\leq l_2(x)\,|u_2-u_1|+l_3(x)\,|v_2-v_1|,
\end{equation}
where
$$
0<x<l,\, u,v, u_i, v_i \in\mathbb{R}, i=1,2,\, \sigma_1(x) \in L_2(0,l),\,\sigma_i(x),\,l_i(x)\in L_\infty(0,l),\, i= 2, 3,
$$
$$
 \sigma_1(x)\geq \text{const}>0, \, \sigma_i(x)\geq 0,\, l_i(x)>0, i=2,3,
$$
are fulfilled.

We impose one more restriction on the beam length $l$ and the parameters $\alpha$ and $\sigma_2(x)$, $\sigma_3(x)$ from the conditions \eqref{eq:bpp-03} and \eqref{eq:bpp-04},\eqref{eq:bpp-05} in the form
\begin{equation}\label{eq:bpp-06}
  \omega =\alpha+\left(\frac{\pi}{l}\right)^2 -\frac{l}{\pi}\left(\frac{2}{\pi}\|\sigma_2(x)\|_\infty+
\|\sigma_3(x)\|_\infty\right) >0.
\end{equation}
Let  us assume that there exists a solution of the problem \eqref{eq:bpp-01},\eqref{eq:bpp-02} and $u\in W_0^{2,2}(0, l)$ [2].
\vskip5mm
\noindent
\section{ The Method}
\vskip 2mm
We will need the Green function for the problem
\begin{equation}\label{eq:bpp-07}
  \begin{gathered}
    v^{\prime\prime\prime\prime}(x)-av^{\prime\prime}(x)=\psi(x),\; \\
    0<x<l, \quad a=\text{const}>0, \\
    v(0)=v(l)=0, \quad v^{\prime\prime}(0)=v^{\prime\prime}(l)=0.
  \end{gathered}
\end{equation}

In order to obtain this function, we split problem \eqref{eq:bpp-07} into two problems
$$
w^{\prime\prime}(x)-aw(x)=\psi(x),
$$
$$
w(0)=w(l)=0,
$$
and
$$
v^{\prime\prime}(x)=w(x),
$$
$$
v(0)=v(l)=0.
$$

Calculations convince us that
$$
w(x)=-\frac{1}{\sqrt{a}\sinh(\sqrt{a}l)}
\Big(
\int_0^x \cosh{\left(\sqrt{a}(x-l)\right)}\,
\cosh{(\sqrt{a}\xi)}\psi(\xi)d\xi
+
$$
$$
+\int_x^l \cosh{(\sqrt{a}x)}
\cosh{(\sqrt{a}(\xi-l))}\psi(\xi)d\xi\Big),
$$
$$
v(x)=\frac{1}{l}\left(
\int_0^x(x-l)\xi w(\xi)d\xi +\int_x^lx(\xi-l)w(\xi)d\xi
\right).
$$

Substituting the first of these formulas into the second and performing integration by parts, we obtain
\begin{equation*}
  \begin{aligned}
    v(x)&=\frac{1}{a} \Big( \int_0^x  \big(k_1(l-x)\xi + k_2  \sinh{(\sqrt{a}(x-l))}\sinh{(\sqrt{a}\xi)} \psi(\xi)\big)\,d\xi + \\
    &+\int_x^l  \big(k_1 x(l-\xi)+k_2  \sinh{(\sqrt{a}x)}\sinh{(\sqrt{a}(\xi-l))} \psi(\xi)\big)\big)\, d\xi  \Big),
  \end{aligned}
\end{equation*}
$$
k_1=
\frac{1}{l},\quad k_2=
\frac{1}{\sqrt{a} \sinh{(\sqrt{a}l)}}.
$$

 The application of \eqref{eq:bpp-07} to problem \eqref{eq:bpp-01}, \eqref{eq:bpp-02} makes it possible to replace the latter problem by the integral equation
 \begin{equation}\label{eq:bpp-08}
   u(x)=\int_0^l G(x,\xi)f(\xi,u(\xi),u^\prime(\xi))\,d\xi,\quad 0<x<l,
 \end{equation}
 where
 $$
 G(x,\xi)=\frac{1}{\tau}                   \begin{cases}  \displaystyle
  \frac{1}{l}(x-l)\xi+\frac{1}{\sqrt{\tau} \sinh{(\sqrt{\tau}\, l)}} \sinh{(\sqrt{\tau}(x-l))} \sinh{(\sqrt{\tau}\xi)}, \quad 0<\xi\leq x<l,\\
  \displaystyle \frac{1}{l}x(\xi-l)+
   \frac{1}{\sqrt{\tau} \sinh{(\sqrt{\tau}\, l)}} \sinh{(\sqrt{\tau} x)} \sinh{(\sqrt{\tau}(\xi-l))}, \quad 0<x\leq \xi<l,
                                            \end{cases}
 $$
   $$
              \tau=m\left(\int_0^l{u^\prime}^2(x)\,dx\right).
   $$

   The equation \eqref{eq:bpp-08} is solved by the method of the Picard iterations. After choosing a function $u_0(x),\; 0\leq x\leq l$, which together with its second derivative vanish for $x=0$ and $x=l$, we find subsequent approximations by the formula
   \begin{equation}\label{eq:bpp-09}
     u_{k+1}(x)=\int_0^l G_k(x,\xi)f(\xi,u_k(\xi),u_k^\prime(\xi))\,d\xi,\quad 0<x<l,\; k=0,1, \cdots
   \end{equation}
 where
$$
 G_k(x,\xi)=
 $$
 $$
 =\frac{1}{\tau_k}                   \begin{cases}  \displaystyle
  \frac{1}{l}(x-l)\xi+\frac{1}{\sqrt{\tau_k} \sinh{(\sqrt{\tau_k}\, l)}} \sinh{(\sqrt{\tau_k}(x-l))} \sinh{(\sqrt{\tau_k}\xi)}, \quad 0<\xi\leq x<l,\\
  \displaystyle \frac{1}{l}x(\xi-l)+
   \frac{1}{\sqrt{\tau_k} \sinh{(\sqrt{\tau_k}\, l)}} \sinh{(\sqrt{\tau_k} x)} \sinh{(\sqrt{\tau_k}(\xi-l))}, \quad 0<x\leq \xi<l,
                                            \end{cases}
 $$
$$
              \tau_k=m\left(\int_0^l{u_k^\prime}^2(x)\,dx\right),
   $$
and $u_k(x)$ is the $\textit k\,$th approximation of the solution of equation \eqref{eq:bpp-08}.

\vskip5mm
\noindent
\section{ The Equation for the Method Error}
\vskip 2mm

Our aim is to estimate the error of the method, by which we understand the difference between the approximate and exact solutions
\begin{equation}\label{eq:bpp-10}
  \delta u_k(x)=u_k(x)-u(x), \quad k=0,1, \cdots .
\end{equation}
For this, it is advisable to use not formula \eqref{eq:bpp-09}, but the system of equalities
\begin{equation}\label{eq:bpp-11}
  u_{k+1}^{\prime\prime\prime\prime}(x)-m\left(\int_0^l {u_k^\prime}^2(x)dx\right)u_{k+1}^{\prime\prime}(x)=f(x,u_k(x),u_k^\prime(x)),
\end{equation}
\begin{equation}\label{eq:bpp-12}
  u_k(0)=u_k(l)=0,\; u_k^{\prime\prime}(0)=u_k^{\prime\prime}(l)=0,
\end{equation}
which follows from \eqref{eq:bpp-09}.

If we subtract the respective equalities in \eqref{eq:bpp-01} and \eqref{eq:bpp-02} from \eqref{eq:bpp-11} and \eqref{eq:bpp-12}, then we get
\begin{equation}\label{eq:bpp-13}
  \begin{aligned}
    &\delta u_k^{\prime\prime\prime\prime}(x)-\frac{1}{2}\Bigg( \bigg[m\left(\int_0^l{u^\prime}^{\,2}_{k-1}(x)\,dx\right)+ m\left(\int_0^l{u^\prime }^{\,2}(x)\,dx\right)\bigg]\delta u_k^{\prime \prime}(x)+ \\
    &+\left [m\left(\int_0^l{u^\prime}^{\,2}_{k-1}(x)\,dx\right)-m\left(\int_0^l{u^\prime}^{\,2}(x)\,dx\right)\right](u^{\prime \prime}_k(x)+u^{\prime\prime}(x)) \Bigg)= \\
    &=f(x,u_{k-1}(x),u^\prime_{k-1}(x))-f(x,u(x),u^\prime(x)),
  \end{aligned}
\end{equation}
\begin{equation}\label{eq:bpp-14}
  \delta u_k(0)=\delta u_k(l)=0, \quad \delta u^{\prime\prime}(0)=\delta u^{\prime\prime}(l)=0, \; k=1,2, \cdots .
\end{equation}

We will come back to \eqref{eq:bpp-13},\eqref{eq:bpp-14} to estimate the error of method \eqref{eq:bpp-09}. In meantime we have to derive several a priori estimates.

\vskip5mm
\noindent
\section{ Auxiliary Inequalities}
\vskip 2mm
Let
\begin{equation}\label{eq:bpp-15}
  \|u(x)\|_p=\left(\int_0^l \left(\frac{d^pu}{dx^p}(x)\right)^2 dx\right)^{1/2},\;\; p=0,1,2, \quad \|u(x)\|=\|u(x)\|_0.
\end{equation}
The symbol $(\cdot,\cdot)$ is understood as a scalar product in $L_2(0,l).$

\begin{lemma}
  The following estimates are true
  \begin{equation}\label{eq:bpp-16}
    \|u(x)\|\leq \frac{l}{\pi} \,\|u(x)\|_1,\quad \|u(x)\|_1\leq \frac{l}{\pi} \,\|u(x)\|_2,
  \end{equation}
  respectively for $u(x)\in W_0^{1,2}(0,l)$ and $u(x)\in W^{2,2}(0,l)\cap W_0^{1,2}(0,l)$.
  \begin{proof}
    Indeed, the first estimate of \eqref{eq:bpp-16} is Friedrich's inequality (see, e.g. [7], p. 192).
    Applying this inequality and taking into account that
    $$
    \|u(x)\|_1^2=u(x)u'(x)|_0^l-(u(x),u''(x))=-(u(x),u''(x))\leq \|u(x)\|\,\|u(x)\|_2
    $$
    we get the second inequality of \eqref{eq:bpp-16}.
  \end{proof}
\end{lemma}
\begin{lemma}
  The inequality
  \begin{equation}\label{eq:bpp-17}
    \|f(x,u(x),u^\prime(x))\|\leq \|\sigma_1(x)\|+\left(\frac{l}{\pi}\|\sigma_2(x)\|_\infty+\|\sigma_3(x)\|_\infty\right)\|u(x)\|_1
  \end{equation}
  is fulfilled for $u(x)\in W_0^{1,2}(0,l).$
  \begin{proof}
    By \eqref{eq:bpp-04} we write
    $$
    \|f(x,u(x),u^\prime(x))\|\leq \|\sigma_1(x)\|+\|\sigma_2(x)\|_\infty\,\|u(x)\|+\|\sigma_3(x)\|_\infty\,\|u^\prime (x)\|
    $$
    Recall also \eqref{eq:bpp-16}. The result is \eqref{eq:bpp-17}.
  \end{proof}
\end{lemma}
\begin{lemma}
  For the solution of problem \eqref{eq:bpp-01},\eqref{eq:bpp-02} we have the inequality
  \begin{equation}\label{eq:bpp-18}
    \|u(x)\|_1\leq c_1,
  \end{equation}
  where
  \begin{equation}\label{eq:bpp-19}
    c_1=\frac{l}{\omega\pi}\|\sigma_1(x)\|.
  \end{equation}
  \begin{proof}
    We multiply  equation \eqref{eq:bpp-01} by $u(x)$ and integrate the resulting equaliti with respect to $x$ from 0 to $l$. Using \eqref{eq:bpp-02}, we get
    $$
    \|u(x)\|^2_2+m(\|u(x)\|^2_1)\|\|u(x)\|_1^2=(f(x,u(x),u^\prime(x)),u(x)).
    $$
    By \eqref{eq:bpp-16} and \eqref{eq:bpp-03} we obtain
    $$
    \left(\alpha+\left(\frac{\pi}{l}\right)^2\right)       \|u(x)\|^2_1\leq \frac{l}{\pi}\|f(x,u(x),u^\prime(x))\|\,\|u(x)\|_1.
    $$
    Therefore by \eqref{eq:bpp-17},
    $$
    \left(
    \alpha+\left(\frac{\pi}{l}\right)^2-
    \left(\frac{l}{\pi}\right)^2\|\sigma_2(x)\|_\infty-\frac{l}{\pi}\|\sigma_3(x)\|_\infty \right)
    \|u(x)\|_1\leq \frac{l}{\pi}\|\sigma_1(x)\|.
    $$
    From this relation and \eqref{eq:bpp-06} follows \eqref{eq:bpp-18}.
  \end{proof}
\end{lemma}
\begin{lemma}
  Suppose where given some numbers $v_k\geq 0, k=0,1,\cdots ,$ for which the inequality
  \begin{equation}\label{eq:bpp-20}
    v_k\leq av_{k-1}+b,\quad k=1,2,\cdots ,
  \end{equation}
  where $0\leq a<1,\; b>0,$ holds. Then we have the following uniform estimate with respect to the index $k$
  \begin{equation}\label{eq:bpp-21}
    v_k \leq \frac{b}{1-a}+a\max \left(0,v_0-\frac{b}{1-a}\right),\; k=1,2,\cdots .
  \end{equation}
  \begin{proof}
    By virtue of \eqref{eq:bpp-20}, by the method of mathematical induction we have $v_k\leq a^kv_0+ (a^{k-1}+a^{k-2}+\cdots +1)b, k=1,2,\cdots $, which implies
    \begin{equation}\label{eq:bpp-22}
      v_k\leq a^kv_0+\frac{1-a^k}{1-a}b=\frac{b}{1-a}+a^k \left(v_0-\frac{b}{1-a}\right).
    \end{equation}
    Let us denote  $\displaystyle\nu_k=a^k\left(v_0-\frac{b}{1-a}\right)$ and consider two cases $\displaystyle v_0\leq \frac{b}{1-a}$ and
    $ \displaystyle  v_0 > \frac{b}{1-a}$. In the first case  $\nu_k\leq 0$ and by virtue of \eqref{eq:bpp-22} $ \displaystyle v_k \leq \frac{b}{1-a}, k=1,2,\cdots .$ In the second case $\displaystyle \nu_k >0, \; \max \nu_k=\nu_1=a(v_0-\frac{b}{1-a}),$  which, by virtue of \eqref{eq:bpp-22} yields $\displaystyle v_k=\frac{b}{1-a}+a(v_0-\frac{b}{1-a}), k=1,2, \cdots .$  From this conclusions the validity of estimate \eqref{eq:bpp-21} follows.
  \end{proof}
\end{lemma}
\begin{lemma}
  Approximations of iteration method \eqref{eq:bpp-09} satisfy the inequality
  \begin{equation}\label{eq:bpp-23}
    \|u_k(x)\|_1\leq c_2, \; k=1,2,\cdots ,
  \end{equation}
  where
  \begin{equation}\label{eq:bpp-24}
    c_2= \begin{cases} c_1, & if \;\|\sigma_2(x)\|_\infty+\|\sigma_3(x)\|_\infty=0,
  \\
      c_1+c_0 \max (0, \|u_0(x)\|_1-c_1),
      & if \;\|\sigma_2(x)\|_\infty+\|\sigma_3(x)\|_\infty\neq 0,
      \end{cases}
  \end{equation}
  \begin{equation*}
    c_0=
  \left(1+\omega
                 \left(\left(\frac{l}{\pi}\right)^2\|\sigma_2(x)\|_\infty
                                  +\frac{l}{\pi}\|\sigma_3(x)\|_\infty\right)^{-1}
  \right)^{-1}.
  \end{equation*}
\begin{proof}
    Replace $k$ by the index $k-1$ in equation \eqref{eq:bpp-11}, multiply the resulting relation by $u_k(x)$ and integrate over $x$ from 0 to $l$. Taking  \eqref{eq:bpp-12} into account, we get
    $$
    \|u_k(x)
    \|^2_2+m(\|u_{k-1}(x)\|^2_1)\,\|u_k(x)\|^2_1=\left(f(x,u_{k-1}(x), u^\prime_{k-1}(x)),\,u_k(x)\right),\; k=1,2,\cdots.
    $$
    Applying \eqref{eq:bpp-03} and \eqref{eq:bpp-15}, we have
    $$
    \left(\alpha+\left( \frac{\pi}{l}\right)^2\right)
    \|u_k(x)\|^2_1 \leq \frac{l}{\pi}\|f(x,u_{k-1}(x),u^\prime_{k-1}(x))\|\,\|u_k(x)\|_1,
    $$
    which implies
    $$
    \left(\alpha+\left( \frac{\pi}{l}\right)^2\right)
    \|u_k (x)\|_1 \leq \frac{l}{\pi}\|f(x,u_{k-1}(x)\,u^\prime_{k-1}(x))\|.
    $$
    Hence, using \eqref{eq:bpp-17}, we conclude that
    $$
    \|u_k(x)\|_1\leq \frac{1}{\alpha+(\frac{\pi}{l})^2}\frac{l}{\pi}
    \left(\|\sigma_1(x)\|+\left(\|\sigma_2(x)\|_\infty\frac{l}{\pi}+\|\sigma_3(x)\|_\infty\right)\|u_{k-1}(x)\|_1\right).
    $$
    This relation is an inequality of type \eqref{eq:bpp-20}, where
    $$
    v_k=\|u_k(x)\|_1,\quad a=\frac{1}{\alpha+(\frac{\pi}{l})^2}\frac {l}{\pi}\big(\|\sigma_2(x)\|_\infty\frac{l}{\pi}+\|\sigma_3(x)\|_\infty\big),\quad
    b=\frac{1}{\alpha+(\frac{\pi}{l})^2}\frac {l}{\pi}\|\sigma_1(x)\|.
    $$
    Let us apply \eqref{eq:bpp-06}, \eqref{eq:bpp-19} to these formulas and carry out some calculations. As a result, for  $\|\sigma_2(x)\|_\infty+\|\sigma_3(x)\|_\infty=0$ we obtain $a=0$ and $\displaystyle \frac{b}{1-a}=c_1$, while for  $\|\sigma_2(x)\|_\infty+\|\sigma_3(x)\|_\infty\neq0$  we have $a=c_0$  and $\displaystyle \frac{b}{1-a}=c_1$.  By considering these two cases with estimate \eqref{eq:bpp-21} we get convinced that \eqref{eq:bpp-23} is valid.
  \end{proof}
\end{lemma}
By Lemma 3 and Lemma 5 it will be natural to require that the initial approximatio  $u_0(x)$ in \eqref{eq:bpp-09} satisfy the condition
\begin{equation}\label{eq:bpp-25}
  \|u_0(x)\|_1\leq c_1.
\end{equation}
Then, by virtue of \eqref{eq:bpp-24} and \eqref{eq:bpp-23}, we have  $\|u_k(x)\|_1\leq c_1$, which, with \eqref{eq:bpp-19} taken into account, implies
\begin{equation}\label{eq:bpp-26}
  \|u_k(x)\|_1\leq \frac{l}{\omega \pi}\|\sigma_1(x)\|, \quad k=0,1,\cdots .
\end{equation}

\vskip5mm
\noindent
\section{ Convergence of the Method}
\vskip 2mm

Multiplying \eqref{eq:bpp-13} by $\delta u_k(x)$, integrating the resulting equality  with respect to $x$ from 0 to $l$ and using \eqref{eq:bpp-14}, we come to the relation
$$
\|\delta u_k(x)\|^2_2 +\frac{1}{2}\Big(\left( m\left(\|u_{k-1}(x)\|^2_1\right) +m\left(\|u(x)\|^2_1
\right)\right)\|\delta u_k(x)\|^2_1+
$$
$$
+\left(m\left(\|u_{k-1}(x)\|^2_1\right)-m\left(\|u(x)\|^2_1\right)\right)\left(u^\prime_k(x)+u^\prime(x), \delta u^\prime_k(x)\right)\Big)=
$$
$$
= \left(f(x,u_{k-1}(x), u^\prime_{k-1}(x))-f(x,u(x),u^\prime(x)), \delta u_k(x)\right), \quad k=1,2, \cdots .
$$
Applying \eqref{eq:bpp-03}-\eqref{eq:bpp-05} and \eqref{eq:bpp-16} we first obtain
$$
\|\delta u_k(x)\|^2_2+ \alpha \|\delta u_k(x)\|^2_1\leq
 \frac{1}{2}l_1 \prod_{p=0}^1\left|\left(u^\prime_{k-p}(x)+u^\prime, \delta u^\prime_{k-p}(x)\right)\right|+
$$
$$
+\Big(\|l_2(x)\|_\infty\,\|\delta u_{k-1}(x)\|+\|l_3(x)\|_\infty\,\|\delta u^\prime_{k-1}(x)\|\Big)\|\delta u_k(x)\|\leq
$$
$$
\leq \frac{1}{2}l_1 \prod_{p=0}^1 \Big(\|u_{k-p}(x)\|_1+\|u(x)\|_1\Big) \| \delta u_{k-p}(x)\|_1+
\frac{l}{\pi}\Big( \frac{l}{\pi}\|l_2(x)\|_\infty+ \|l_3(x)\|_\infty \Big) \prod_{p=0}^1\|\delta u_{k-p}(x)\|_1,
$$
and after that, by virtue of \eqref{eq:bpp-18} and \eqref{eq:bpp-26} we have
$$
\|\delta u_k(x)\|_1\leq
 \Big(\alpha+\big(\frac{\pi}{l}\big)^2\Big)^{-1}
 \Big(\frac{1}{2}l_1\prod_{p=0}^1\big(\|u_{k-p}(x)\|_1+\|u\|_1(x)\big)+
$$
$$
+\big(\frac{l}{\pi}\big)^2\|l_2(x)\|_\infty+\frac{l}{\pi}\|l_3(x)\|_\infty\Big)\|\delta u_{k-1}\|_1\leq q\|\delta u_{k-1}\|_1,\quad
k=1,2, \cdots ,
$$
where
$$
q= \left(\alpha +\big(\frac{\pi}{l}\big)^2\right)^{-1}
\Big(2c_1^2l_1+\|l_2(x)\|_\infty\big(\frac{l}{\pi}\big)^2+
\|l_3(x)\|_\infty \frac{l}{\pi}\Big).
$$
Taking \eqref{eq:bpp-10},\eqref{eq:bpp-19} and \eqref{eq:bpp-16} into consideration we come to the following result
\begin{thm}
  Let assumptions \eqref{eq:bpp-03}-\eqref{eq:bpp-06} and \eqref{eq:bpp-25} are fulfilled. Suppose besides
  $$
  q=\frac{1}{\alpha+\big(\frac{\pi}{l}\big)^2}\left(\frac{l}{\pi}\right)^2 \Big(2l_1\big(\frac{\|\sigma_1(x)\|}{\omega}\big)^2+
  \|l_2(x)\|_\infty+\frac{\pi}{l}\|l_3(x)\|_\infty\Big)<1.
  $$
  Then the approximations of the iteration method \eqref{eq:bpp-09} converge to exact solution of problem \eqref{eq:bpp-01},\eqref{eq:bpp-02} and for the error the folloving estimate
  $$
  \|u_k(x)-u(x)\|_p \leq \Big(\frac{l}{\pi}\Big)^{1-p}q^k\|u_0(x)-u(x)\|_1,\; k=1,2, \cdots ,\; p=0,1 ,
  $$
  is true.
\end{thm}

\vskip5mm
\noindent
\section{ Numerical Experiment}
\vskip 2mm

The theoretical results about the convergence of approximations of iteration method \eqref{eq:bpp-09} to exact solution ${u}{\left( x \right)}$ of problem \eqref{eq:bpp-01}, \eqref{eq:bpp-02} is confirmed. For illustration, the results of numerical computations of one of the test problems are given below.

We consider a special case, where $m\left( z \right)={{m}_{0}}+{{m}_{1}}\cdot z,\,\,{{m}_{0}},{{m}_{1}}>0,\,\,$ ${{m}_{0}}=1,\,\,{{m}_{1}}=\frac{1}{2},$ the beam length $l=1,$ exact solution $u\left( x \right)=x\left( x-1 \right)\left( {{x}^{2}}-x-1 \right),$ i. e. $u\left( x \right)={{x}^{4}}-2{{x}^{3}}+x,$ the right-hand side
\begin{align*}
  f\left( x,u\left( x \right),{u}'\left( x \right) \right)&=\frac{1}{35} \big (43.5{{{{u}'}}^{2}}\left( x \right)-348{{x}^{3}}{u}'\left( x \right)-1566u\left( x \right)+ \\
  & +696{{x}^{6}}-3132{{x}^{3}}+2088x+796.5 \big )\,.
\end{align*}

We carried out five, seven and nine iterations. To compute the integrals on $\left[ 0,1 \right]$ we divided the interval into $n=10, 20$ parts $\left( h=0.1,\,\,\,0.05,\,\,\,\text{respectively} \right)$ and used the square formula of trapezoid. The error in the k--iteration is defined as
\begin{equation*}
  error\,k=\,\underset{i=0,1,...,n}{\mathop{\max }}\,\,\left\{ \left| {{u}_{k}}\left( {{x}_{i}} \right)-u\left( {{x}_{i}} \right) \right| \right\},\,\,{x}_{i}={i}{h},\,\,k=1,2,...,9.
\end{equation*}
Numerical values for the errors are calculated (see Table 1).

The function ${{u}_{0}}\left( x \right)=0$ is taken as the initial approximation. In case of five, seven and nine iterations for $n=10,\,20$ the exact and approximate solutions are graphically illustrated (Figs. 1-6).
\begin{remark}
  In the figures the green line color denotes the exact solution graph, yellow is the first approximation, red -- the second, blue -- the third, pink -- the fourth, golden -- the fifth, brown -- the sixth, purple -- the seventh, orange -- the eighth and black -- the ninth.
\end{remark}
\begin{table}[h!]
  \centering
  \begin{tabular}{|c|c|c|c|c|c|c|c|}
  \hline
  $n \backslash error$ & $error\,1$ & $error\,2$ & $error\,3$ & $error\,4$ & $error\,5$ & $error\,7$ & $error\,9$ \\
  \hline
  $n=10$ & $0.43203$ & $0.16734$ & $0.06405$ & $0.02473$ & $0.00953$ & $0.00142$ & $0.00021$ \\
  \hline
  $n=20$ & $0.43328$ & $0.16715$ & $0.06365$ & $0.02446$ & $0.00938$ & $0.00138$ & $0.00020$ \\
  \hline
\end{tabular}
  \caption{}\label{table:1}
\end{table}\newpage

\begin{figure}[h!]
  \centering
  \begin{minipage}[b]{0.495\textwidth}
    \includegraphics[width=\textwidth]{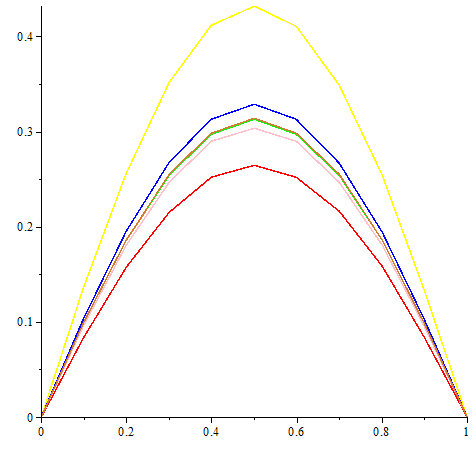}
    \caption{$Iteration=5,\,n=10$}
  \end{minipage}
  \hfill
  \begin{minipage}[b]{0.495\textwidth}
    \includegraphics[width=\textwidth]{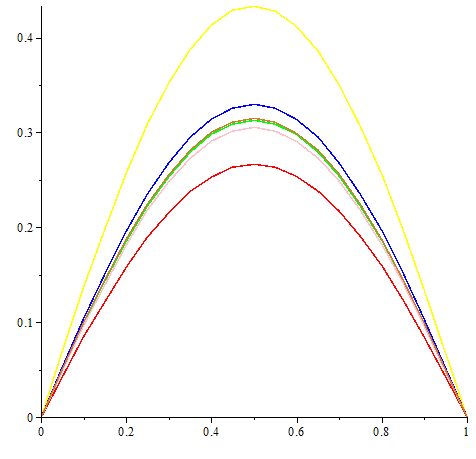}
    \caption{$Iteration=5,\,n=20$}
  \end{minipage}
\end{figure}
\begin{figure}[h!]
  \centering
  \begin{minipage}[b]{0.495\textwidth}
    \includegraphics[width=\textwidth]{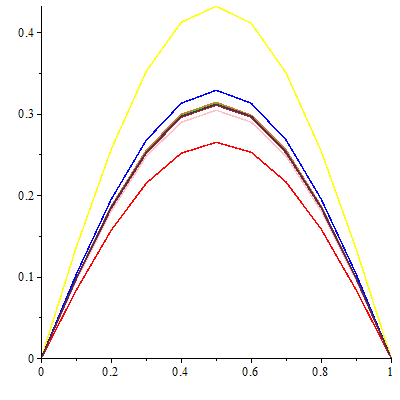}
    \caption{$Iteration=7,\,n=10$}
  \end{minipage}
  \hfill
  \begin{minipage}[b]{0.495\textwidth}
    \includegraphics[width=\textwidth]{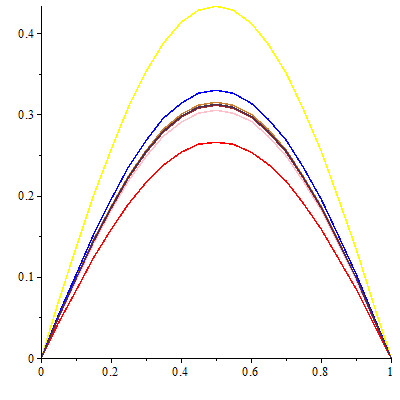}
    \caption{$Iteration=7,\,n=20$}
  \end{minipage}
\end{figure}\newpage
\begin{figure}[h!]
  \centering
  \begin{minipage}[b]{0.495\textwidth}
    \includegraphics[width=\textwidth]{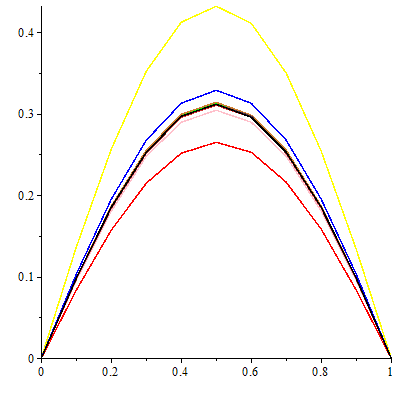}
    \caption{$Iteration=9,\,n=10$}
  \end{minipage}
  \hfill
  \begin{minipage}[b]{0.495\textwidth}
    \includegraphics[width=\textwidth]{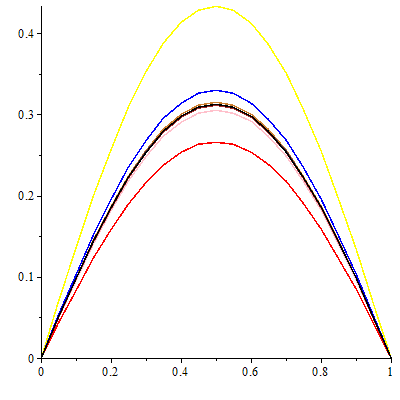}
    \caption{$Iteration=9,\,n=20$}
  \end{minipage}
\end{figure}
The numerical experiments clearly show the convergence of iteration approximate solutions to the exact solution of the problem. The error decreases with the growth of the parameters $n$ and $k$.

\vskip 8mm
\textbf{References}
\vskip 3mm

[1] C. Bernardi and M.I.M. Copetti, Finite element discretization of a thermoelastic beam. \textit{Archive Ouverte HAL-UPMC}, 29/05/2013, 23pp.
\vskip 3mm

[2] S. Fu\v{c}ik and A. Kufner, Nonlinear differential equations. Studies in Applied Mechanics, 2.\textit{ Elsevier Scientific Publishing Company, Amsterdam-Oxford-New York,} 1980.
\vskip 3mm

[3] G. Kirchhoff, Vorlesungen \"{u}ber mathematische physik, I. Mechanik. \textit{Teubner, Leipzig,} 1876.
\vskip 3mm

[4] T.F. Ma, Positive solutions for a nonlocal fourth order equation of Kirchhoff type. \textit{Discrete Contin. Dyn. Syst.} 2007, 694--703.
\vskip 3mm

[5] J. Peradze, A numerical algorithm for a Kirchhoff-type nonlinear static beam. \textit{J.Appl. Math. }\textbf{2009}, Art.ID 818269, 12pp.
\vskip 3mm

[6] J. Peradze, On an iteration method of finding a solution of a nonlinear equilibrium problem for the Timoshenko plate. \textit{ZAMM Z. Angew. Math. Mech.}\textbf{91} (2011), no. 12, 993 --1001.
\vskip 3mm

[7] K.Rektorys, Variational methods in mathematics, science and engineering.\textit{ Springer Science $\&$ Business Media,} 2012.
\vskip 3mm

[8] H. Temimi, A.R. Ansari and A.M.
Siddiqui, An approximate solution for the static beam problem and nonlinear integro-differential equations.\textit{Comput. Math. Appl.}\textbf{62} (2011). no. 8, 3132--3139.
\vskip 3mm

[9] S.Y. Tsai, Numerical computation for nonlinear beam problems.\textit{M.S. thesis, National Sun Yat-Sen University, Kaohsiung, Taiwan,} 2005.
\vskip 3mm

[10] F. Wang and Y. An, Existence and multiplicity of solutions for a fourth-order elliptic equation. \textit{Bound. Value Probl.} \textbf{2012}, 2012:6, 9 pp.
\vskip 3mm

[11] S.Woinowski-Krieger, The effect of an axial force on the vibration of hinged bars. \textit{J. Appl. Mech.}\textbf{17} (1950), 35--36.

\end{document}